\documentclass[a4paper,11pt]{article}
\usepackage{amsfonts,amssymb,latexsym,amsmath}
\begin{document}
\Large
\newcounter{num}[section]
\setcounter{num}{0}
\renewcommand{\thenum}{\arabic{section}.\arabic{num}}

\newcommand{\Num}{\refstepcounter{num}%
\textbf{\arabic{section}.\arabic{num}}}

\newcommand{\Lemma}{\textbf{Lemma~}}
\newcommand{\Prop}{\textbf{Proposition~}}
\newcommand{\Cor}{\textbf{Corollary~}}
\newcommand{\Theorem}{\textbf{Theorem~}}
\newcommand{\Proof}{\textbf{Proof}}

\newcommand{\De}{\Delta}
\newcommand{\al}{\alpha}
\newcommand{\la}{\lambda}
\newcommand{\vpi}{\varpi}
\newcommand{\veps}{\varepsilon}
\newcommand{\gog}{{\mathfrak g}}
\newcommand{\nog}{{\mathfrak n}}
\newcommand{\mog}{{\mathfrak m}}
\newcommand{\hog}{{\mathfrak h}}
\newcommand{\Aog}{{\mathfrak A}}

\newcommand{\AC}{{\cal A}}
\newcommand{\BC}{{\cal B}}
\newcommand{\FC}{{\cal F}}
\newcommand{\BAC}{{{\cal B}{\cal A}}}
\newcommand{\BFC}{{{\cal B}{\cal F}}}

\newcommand{\bog}{{\mathfrak b}}

\newcommand{\eps}{\varepsilon}
\newcommand{\ad}{{\mathrm{ad}}}

\newcommand{\Ad}{{\mathrm{Ad}}}
\newcommand{\rank}{{\mathrm{rank}}}

\newcommand{\Zb}{{\Bbb Z}}
\newcommand{\Nb}{{\Bbb N}}
\newcommand{\Ab}{{\Bbb A}}

\renewcommand{\leq}{\leqslant}
\renewcommand{\geq}{\geqslant}
\newcommand{\xt}{\tilde{x}}
\newcommand{\gt}{\tilde{g}}

\date{}
\title{Reduction of spherical functions}
\author{A.N.Panov
\footnote {The paper is supported by RFBR grants 08-01-00151-a,
09-01-00058-a and by ADTP grant 3341}} \maketitle
 \begin{center}
{\it\large Department of Mathematical Sciences, Samara State
University,\\
443011, Samara, ul. Akademika Pavlova 1, Russia\\
E-mail: apanov@list.ru }
\end{center}

\begin{abstract}  Using  reduction of spherical functions we obtain
generators of the algebra and the field of invariants for the
coadjoint representation of Borel and maximal nilpotent subalgebras
of simple Lie algebras.
\end{abstract}

\section{Introduction}

In this paper we study the rings and fields of invariants of the
coadjoint representation of Borel and maximal nilpotent subalgebras
of simple Lie algebras.

In the original paper on the orbit method ~\cite{K1}, the generators
of the ring of polynomial invariants of the coadjoint representation
of unitriangular Lie algebra (i.e. maximal nilpotent subalgebra in
$A_n$) were discovered. It was proved that the algebra of invariants
is freely generated by the system of corner minors (see also
~\cite{K-Orb}).

There are no nontrivial polynomial invariants for the coadjoint
representation of  Borel subalgebra (see theorem \ref{B-inv}). But
the field of invariants may be nontrivial, for instance for $A_n$.
The system of generators of the field of invariants for the Borel
subalgebra of $A_n$ was discovered in the papers ~\cite{Deift, Arch}
(see also ~\cite[4.8]{P}). In the paper ~\cite{Trof}, solving the
problem of integration of Eiler equation,  the generators for the
Borel subalgebras  of simple Lie algebras $B_n$,~ $C_n$, $D_n$ and
$G_2$ were fond out. For the classical Lie algebras the statement
was formulated in terms of the minors of characteristic matrix.

In the general case of an arbitrary simple Lie algebra, one can
obtain the description of generators as a result of some induction
procedure (see proposition \ref{LocZ} and ~\cite{Gekht}). However,
it will be desirable to obtain more exact description, as for
classical Lie algebras. But it is not clear in what terms to
formulate this exact description. How to determine what is a minor
 and a characteristic matrix for an arbitrary  simple Lie algebra?

We offer the  new approach in this paper, that is gives a
possibility to obtain generators of the ring and the field of
invariants for the Borel and  maximal nilpotent subalgebras.

The approach is based on the reduction of spherical functions. The
mail results of the paper are formulated in the theorems  2.12, 3.1,
3.2, 3.6, 3.7.

\section{Invariants of the coadjoint representation of maximal nilpotent subalgebras}

Let $\gog$ be a split Lie algebra over a field  $K$ of zero
characteristic with the roots system ~ $\De$. By $G$ we define the
linear subgroup with the Lie algebra  $\gog$. Introduce the
following
notations:\\
 $\hog$ is the standard Cartan subalgebra;\\
  $\nog$ (resp. $\nog_-$) is the nilpotent subalgebra, spanned by root vectors  $e_\al$, ~$\al >0$ (resp. $\al<0$);\\
 $\bog$ (resp. $\bog_-$) is the Borel subalgebra, which is equal to  $\bog =
\hog\oplus\nog$ (resp. $\bog_- = \hog\oplus\nog_-$);\\
 $H$,~ $N$,~
$N_-$, $B$,~ $B_-$  are the subgroups of $G$ corresponding to these
subalgebras;\\
 $\nog^*$ and $\bog^*$ are conjugate subspaces for $\nog$ and
$\bog$;\\
 $\AC$ is the algebra  $K[\nog^*]$ of regular functions on  $\nog^*$;\\
$\BAC$ is the algebra  $K[\bog^*]$ of regular functions on $\bog^*$;\\
$\FC$ is the field $K(\nog^*)$ of rational functions on  $\nog^*$;\\
 $\BFC$ is the field   $K(\bog^*)$ of rational functions on  $\bog^*$.

Recall that the coadjoint representation of the group  $N$ (the same
for any other Lie group) is defined by the formula
$$\Ad^*_g f(x_+) = f(\Ad^{-1}_gx_+),$$ where $f\in \nog^*$, ~$x_+\in \nog$,
~$g\in N$. Using the Killing form  $(\cdot,\cdot)$, we identify
$\nog_-$ (resp. $\bog_-$) with the conjugate space  $\nog^*$ (resp.
$\bog^*$).
  Taking into account identification $\nog_*$ with
$\nog_-$, we have got $$\Ad^*_g(x) = \pi_-(\Ad_g(x)),$$ where $x\in
\nog_-$, ~ $g\in N$ and $\pi_-$ is the natural projection $\gog$ on
$\nog_-$.

 The algebra  $\AC$ is a Poisson algebra with respect to linear Poisson bracket such that  $\{x,y\} = [x,y]$
for any  $x,y\in\nog$. The symplectic leaves  of this Poisson
bracket coincide with the coadjoint orbits of the group  $N$ in
$\nog^*$ ~\cite{K1, K-Orb}.

In this section we shall obtain the generators in the algebra of
invariants $\AC^N$  of the coadjoint representation of the group
$N$.

To any rational function  $F(g)$ on the group  $G$ we correspond the
formal series
\begin{equation} \label{Fo} F(\exp tx) = t^k(F_0(x) + tF_1(x)+ t^2 F_2(x)
+\ldots ) , \end{equation}
 where $ k\in \Zb$ and the coefficients  $F_i(x)$  are rational functions on  $\nog_-$.  If $F$ belongs to the local ring
 of the unit element $e$ (in particular,  $F\in K[G]$), then $k\in
\Zb_+$ and all coefficients are lying in  $K[\nog_-]$. We shall say
that
 $F_0$ is the lowest coefficient in decomposition (\ref{Fo}) and $k$ is the lowest degree.

Denote by  $K[G]^{N\times N}$ the ring of invariants with respect to
left-right action of the group $N\times N$ in  $K[G]$.\\
\\
\Prop\Num\label{NN}. {\it If $F\in K[G]^{N\times N}$, then $F_0\in
\AC^N$.}\\
\\
\Proof. First of all we note that it is sufficient to prove the
statement in the case when  $F$  is a eigenfunction with respect to
the right action of the Borel subgroup $B$. That is when
$F(gb)=\chi(b)F(g)$, where  $g\in G$,~ $b\in B$ and $\chi$ is some
character of the subgroup
 $B$.

 Since  $\gog=\nog_-+\bog$, then $N_-B$  is a Zarisky-open subset in $G$. Any element $a\in N_-B$
is uniquely presented in the form  $a=a_-a_+$, where $a_-\in N_-$
and $a_+\in B$.

For any  $s\in N_-$ there exists an open neighbourhood  of unity
such that  $gs\in N_-B$ for any  $g$ in this neighbourhood and,
therefore,~ $gs=(gs)_-(gs)_+$. Denote
$$\rho_g(s) = (gs)_-.$$ This formula determines the local action of
$G$ on $N_-$, which is called the  dressing action. In particular,
if we put  $g\in N$, then this formula determines the dressing
action of $N$ on  $N_-$.

Since $F\in K[G]^{N\times N}$, then for any  $g\in N$ and $s\in N_-$
we have
\begin{equation}\label{gs}
F(\rho_g(s)) = F((gs)_-) = \frac{F((gs)_-(gs)_+)}{\chi((gs)_+)} =
\frac{F(gs)}{\chi((gs)_+)} = \frac{F(s)}{\chi((gs)_+)}.
\end{equation}

Substitute $s=\exp(tx)$, where $x\in \nog_-$, into the formula
(\ref{gs}).
 Since   $\chi(g) = 1$ for any  $g\in N$, then
$$\chi((g\exp(tx))_+) = 1+t\theta(x,t),$$
for some polynomial  $\theta(x,t)$.

  Denote by  $\eta(t)$ the curve  $\rho_g(\exp(tx))$ in the group $N_-$. The formula (\ref{gs})
implies
\begin{equation} \label{et} F(\eta(t))
=\frac{F(\exp(tx))}{1+t\theta(x,t)}.
\end{equation}

As $\eta(0) = \left.\exp(tx)\right|_{t=0} = e$, then, by the formula
(\ref{et}), we obtain
\begin{equation} \label{eto} F_0(\eta'(0)) =
F_0(x).
\end{equation}

Let us show that  $\eta'(0) = \Ad_g^*(x)$. Since  $g\in N$, then for
small  $t$ the curve $g\exp(tx)$ belongs to the open subset $N_-B$.
From this it follows that
$$ g\exp(tx) = \eta(t)\zeta_1(t),$$
for some curve  $\zeta_1(t)$ in the group $B$. Note that $\zeta_1(0)
= g$. Denoting  $\zeta(t) = \zeta_1(t)g^{-1}$, we have got
\begin{equation}\label{zet}
g\exp(tx)g^{-1} = \eta(t)\zeta(t),
\end{equation}
where  $\zeta(0)=e$. Differentiate  (\ref{zet}) with respect to $t$
at $t=0$:

$$ \Ad_g(x) = \left.\frac{d}{dt} g\exp(tx)g^{-1}\right|_{t=0} =
\eta'(0)\zeta(0) + \eta(0)\zeta'(0) = \eta'(0) + \zeta'(0).$$ As
$\eta'(0)\in \nog_-$ and $\zeta'(0)\in \bog$, then
$$ \eta'(0) = \pi_-(\Ad_g(x)) = \Ad_g^*(x).$$
Substituting into (\ref{eto}), we obtain  $F_0 (\Ad^*_g(x)) =
F_0(x)$ for any  $x\in \nog_-=\nog^*$. We proved that  the
polynomial  $F_0$ is invariant with respect to  $\Ad_g^*$. $\Box$

For any irreducible finite-dimensional  representation  $T$ we
denote by $S_T(g)$ the spherical function
$$S_T(g) = l_0(T_g v_0),$$ where $v_0$ (resp. $l_0$) is the  dominant vector of representation
 $T$ (resp. of conjugate representation of  $T$).\\
 \Cor\Num.~ {\it  For any irreducible finite-dimensional  representation  $T$
the lowest coefficient  $(S_T)_0$ of the expansion (\ref{Fo}) for
$S_T(g)$ belongs to  $\AC^N$}.

Let  $T_1, \ldots, T_n$ be the system of fundamental representations
of the group $G$, their fundamental weights are  $\vpi_1, \ldots,
\vpi_n$, where $n=\rank(\gog)$. Let  $S_1(g),\ldots, S_n(g)$ be the
corresponding spherical functions. For any  $i$ the expansion
(\ref{Fo}) has the form:
\begin{equation}\label{Sx} S_i(\exp tx) = t^k(S_{i0}(x) + tS_{i1}(x)+ t^2 S_{i2}(x)
+\ldots ).\end{equation} By $P_1,\ldots, P_n$ we denote the lowest
coefficients  $S_{10},\ldots, S_{n0}$ of the corresponding decompositions  (\ref{Sx}).\\
\\
  \Cor\Num. {\it The polynomials  $P_1,\ldots, P_n$ belong to $\AC^N$}.

Let $\Gamma$ be the Heisenberg algebra over the field  $K$,
generated by the system  $$x_1,\ldots,
 x_n, y_1,\ldots, y_n, z, $$ satisfying the relations
 $[x_i, y_j]=\delta_{ij} z$ and $[x_i, z] = [y_j, z] =0$.
 Denote by  $V$ the linear subspace, spanned by
 $x_i,~y_j$, where $i,j = \overline{1,n}$. The symmetric algebra
 $\mathrm{Sym}(\Gamma)$ is a Poisson algebra with respect to the linear Poisson bracket, that coincides with commutator on
  $\Gamma$.

 We say that  $\Ab_n$  is a standard Poisson algebra  if it is generated by  $$p_1, \ldots, p_n, q_1, \ldots, q_n$$ satisfying the relations
  $\{p_i, q_j\} = \delta_{ij}$ and
 $\{p_i, p_j\} = \{q_i, q_j\} = 0$. Note that the  localization of  $\mathrm{Sym}(\Gamma)$ with respect to
 $z$ contains the standard Poisson subalgebra   $\Ab_n$,  generated by
 $p_i = x_i$ and $q_j = z^{-1}y_j$.

Recall that a Poisson algebra  $\BC$ is the tensor product of two
its
 Poisson subalgebras
$\BC_1$, ~$\BC_2$, if  $\BC = \BC_1\otimes\BC_2$ as a commutative
algebra  and $\{\BC_1,\BC_2\} =0$.

The localization of $\mathrm{Sym}(\Gamma)$ with respect to
 $z$ is a tensor product  $K[z^{\pm 1}]\otimes \Ab_s$ as a Poisson algebra.\\
  \textbf{Definition}. One says the a linear mapping  $D: \BC\to\BC$
 is a derivation of the  Poisson subalgebra $\BC$, if  $$D(ab) =
 D(a)b + aD(b),$$
 $$D\{a,b\} = \{D(a),b\} + \{a,D(b)\}$$
 for any  $a,b\in \BC$.\\
 \\
  \Lemma \Num\label{diff}. {\it
  For any derivation $D$  of the Poisson algebra
  $\mathrm{Sym}(\Gamma)$, such that $D(V)\subseteq V$ and $D(z) = 0$,  there exists a unique element
   $a_D\in z^{-1} V\mathrm{Sym}(V)$, satisfying  $D(P) =
  \{a_D,P\}$ for any  $P\in \mathrm{Sym}(\Gamma)$}.
  \\
  \\
 \Proof. One can extend the derivation  $D$ to  derivation of  $\Ab_n$ (the standard Poisson subalgebra
 in the localization of  $\mathrm{Sym}(\Gamma)$ with respect to $z$). Any derivation of a standard Poisson subalgebra is
 inner (one can prove this similarly to  ~\cite[4.6.8]{Dix}).
 There exists an element $a\in \Ab_n$, such that  $D(u) = \{a,u\}$ for any
 $u\in \Ab_n$. Since  $D(V)\subseteq V$ and $D(z)=0$, then the element  $a$ can be represented in the form
 $a=z^{-1}b$, where $b\in V\mathrm{Sym}(V)$; this proves the existence of  $a_D$. Easy to see that $a_D$ is uniquely determined by  $D$. $\Box$

Denote by  $e_\alpha$,~ $\al\in \De_+$, the standard basis in
$\nog$. Each vector of the basis (as well as any vector of $\nog$)
is a linear form on  $\nog^*$  and, therefore, it is an element in
$\AC$.

Let $\xi_1$ be the greatest  root in $\De^+$. Denote by  $Z_1$ the
basic vector of weight  $\xi_1$. The Poisson bracket of $\AC$ may be
extended to its localization $\AC(\xi_1)$ with respect to the system
of denominators $\{Z_1^m:~ m\in \Nb\}$.

We say that a root  $\al\in \De^+$ is singular for a root
$\gamma\in \De^+$, if $\gamma-\al\in \De^+$. The subspace
$\Gamma_1$, spanned by $Z_1$ and by all vectors  $e_\al$, where
$\al$  is a singular root for $\xi_1$, is a  Heisenberg algebra.
Denote by  $\De_1^+$ the subset of positive roots, that can be
obtained after deleting in $\De^+$ the root  $\xi_1$ and all
singular for  $\xi_1$ roots. The subspace $\nog_1$, spanned by
$e_\al$, where $\al\in
\De_1^+$, is a Lie subalgebra.\\
\\
 \Lemma\Num\label{xi}.~ {\it The Poisson algebra  $\AC(\xi_1)$ is isomorphic
 to the tensor product  $K[Z_1^\pm]\otimes\Ab_s\otimes \AC_1$
of the commutative Poisson algebra $K[Z_1^\pm] $, some standard
Poisson algebra  $\Ab_s$ and the Poisson algebra  $\AC_1
= K[\nog_1^*]$.}\\
\\
 \Proof. As it was proved above, the localization of the Poisson algebra  $\mathrm{Sym}(\Gamma_1)$ with respect to  $Z_1$ is the tensor product
 $K[Z_1^{\pm 1}]\otimes \Ab_s$.  In our case, the subspace  $V$ of the lemma \ref{diff} is spanned by the root vectors
 $e_\al$, where  $\al$ runs over the set of singular for $\xi_1$
 roots. For any root  $\beta\in\De_1^+$ we consider the derivation
$$D_\beta = \ad_{e_{\beta}}$$ of the Heisenberg algebra $\Gamma_1$. By lemma  \ref{diff},   the exists the element  $a_\beta\in Z_1^{-1}
V_1\mathrm{Sym}(V_1)$ such that  $D_\beta(P) = \{a_\beta, P\}$ for
any $P\in \Gamma_1$. Therefore, the element $\tilde{e}_\beta =
e_\beta - a_\beta$ satisfies
$$\{\tilde{e}_\beta, P\} = 0$$ for any $P\in\Gamma_1$.

The uniqueness of $a_\beta$ implies that the correspondence  $
e_\beta \to \tilde{e}_\beta$   uniquely extends to the embedding of
the Poisson algebra  $\AC_1$ into the Poisson algebra  $\AC(\xi_1)$,
such that its image lies in involution for  $\Gamma_1$. One can see
that  $\AC(\xi_1)$, as a commutative algebra, is generated by
$\Gamma_1$ and the image of  $\AC_1$. This concludes the proof.
$\Box$

In the extended Dynkin diagram the root  $-\xi_1$ connects with one
or two simple roots (the last only for  $A_n$).  After removing
these roots in the system of simple roots  $\Pi$
 of the algebra  $\gog$ we obtain the system of simple roots $\Pi_1$, which is irreducible for all systems of simple roots, besides $B_n$ and $D_n$.
 In the last  cases, the system $\Pi_1$ is the union of  $A_1$ and the corresponding  simple root system of rank  $n-2$ (besides  $D_4$, when this system is a union of
  three $A_1$) .
 The algebra  $\nog_1$ of lemma
 \ref{xi}
is a maximal nilpotent subalgebra in the semisimple  Lie algebra
with the system of simple roots  $\Pi_1$. We choose the greatest
positive root $\xi_2$ for $\nog_1$, if $\Pi_1$ is irreducible, or a
pair of maximal positive roots  $\xi_2>\xi_3$ (resp. triple
$\xi_2>\xi_3>\xi_4$ for $D_4$), if $\Pi_1$ is reducible.
 Continuing this process further, we have got the subsystem of positive roots
 \begin{equation}\label{Xi} \Xi =\{\xi_1, \xi_2, \xi_3, \ldots, \xi_m \}.\end{equation}
 \Prop\Num\label{LocZ}. {\it There exists a system of rational functions   $Z_1,\ldots,
Z_m$ with the system of weights  $\xi_1,\ldots, \xi_m$ with respect
to the coadjoint action of Cartan subalgebra such that \\
1) any $Z_i$ lies in the localization of the algebra  $\AC$ with
respect to the denominator system, generated by  $Z_1,\ldots, Z_{i-1}$ (in particular, ~$Z_1\in\AC$);\\
2) all $Z_i$ are invariant with respect to the coadjoint
representation  of the group $N$ (i.e. all are lying in $\FC^N$);\\
3) the localization  $\AC(\Xi)$ of the algebra $\AC$  with respect
to the denominator system, generated by  $Z_1,\ldots, Z_{m}$, is
isomorphic as a Poisson algebra to the tensor product
$$K[Z^\pm_1,\ldots, Z^\pm_m]\otimes \Ab_s$$ of the commutative Poisson algebra $K[Z^\pm_1,\ldots, Z^\pm_m]$ and some standard  Poisson algebra
 $\Ab_s$.}
\\
\\
\Proof~ is follows from lemma \ref{xi}. $\Box$\\
\\
\Cor\Num \label{zfield}.~ {\it $\FC^N = K(Z_1,\dots,Z_m).$}\\
 \Cor \Num\label{prod}. {\it
 If  $F\in \FC^N$   and $F$ is an eigenfunction for the action of  $H$ in $\FC$, then $F$ is written in the form
\begin{equation}\label{FZ} F = Z_1^{k_1}\cdots
Z_m^{k_m},\end{equation} where $k_i\in\Zb$ for any  $i =
\overline{1,m}$.}
\\
\Proof.~ The system of weights  $\Xi$ is linear independent over $\Zb$. $\Box$\\

This implies that any polynomial  $P_i$ (see (\ref{Sx}))  is written
in the form
\begin{equation} \label{pz} P_i=Z_1^{k_{i1}}\cdots
Z_m^{k_{im}},\end{equation}  where $k_{ij}\in\Zb$.

 The weight of the polynomial  $P_i$ with respect to the action of Cartan subalgebra
$\hog$ is equal to
$$\vpi'_i = (1-w_0)\vpi_i,$$
where $w_0$ is the element of the Weyl group $W$ of greatest length.

Recall that $w_0 = -\mathrm{id}$ for the Lie algebras  $A_n$, $B_n$,
$C_n$, $D_n$ (for even $n$), $G_2$, $F_4$, $E_7$, $E_8$. In all
other cases,  $w_0= - \phi$, there  $\phi$ is the nontrivial
 automorphism of the system of simple roots ~\cite[Tables I-IX]{Bur}: \\
$\phi(\al_i) = \al_{n-i+1}$ in the case $A_n$;\\
$\phi(\al_{n-1}) = \al_n$, ~ $\phi(\al_i) = \al_i $ for $1\leq i\leq
n-2$
in the case $D_n$ (where $n$ is odd);\\
$\phi(\al_1) = \al_6$,~ $\phi(\al_3) = \al_5$, ~ $\phi(\al_2) =
\al_2$,~ $\phi(\al_4) = \al_4$ in the case $E_6$.

The formula (\ref{pz}) implies
\begin{equation}\label{vpk}\vpi'_i= k_{i1}\xi_1 + \ldots +
k_{im}\xi_m.\end{equation}

Further, we find the system $\Xi = \{\xi_1,\ldots,\xi_m\}$ and
obtain the formulas  (\ref{vpk}) for each simple Lie algebra of
types $A_n - E_8$.
In what follows, we use the standard notations of  ~\cite[Tables I-IX]{Bur}.\\
{\textsc Case} $A_n=\mathrm{sl}(n+1,K)$,~ $n\geq 1$. The system
$\Xi$ is consists of  $m=\left[\frac{n+1}{2}\right]$ roots
$$\left\{\begin{array}{l}\xi_1=\veps_1,\\
\xi_2 = \veps_1+\veps_2\\
\ldots \\
\xi_m = \veps_1 +\veps_2+\ldots+\veps_m.\end{array}\right.$$ Thus,
we have got
$$\left\{\begin{array}{l}\vpi_1'= \vpi_n' = \vpi_1+\vpi_n = \xi_1,\\
\vpi_2'= \vpi'_{n-1} = \vpi_2+\vpi_{n-1} = \xi_1 +\xi_2,\\
\ldots\\
\vpi'_m = \vpi_{n-m+1}' = \vpi_m+\vpi_{n-m+1} = \xi_1
+\xi_2+\ldots+\xi_m.\end{array}\right.$$ {\textsc Case}
$B_n=\mathrm{o}(2n+1,K)$, ~$n\geq 2$. Here $m=n$ and
$\vpi'_i=2\vpi_i$ for any  $1\leq i\leq n$. For $n=2$ we  have
$$\left\{\begin{array}{l}\xi_1 = \al_1+2\al_2,\\
\xi_2 = \al_1. \end{array}\right. \quad\quad\quad
\left\{\begin{array}{l}\vpi'_1 = \xi_1+\xi_2,\\
\vpi_2' = \xi_1.
\end{array}\right.
$$

For $n=2l$, ~ $l>1$,  we obtain
$$\left\{\begin{array}{l}\xi_1= \veps_1 + \veps_2,\\
\xi_2 =  \veps_1 - \veps_2,\\
\ldots,\\
\xi_{2l-1}= \veps_{2l-1} +  \veps_{2l},\\
\xi_{2l} =  \veps_{2l-1} -
\veps_{2l}.\end{array}\right.\quad\quad\quad
\left\{\begin{array}{l}\vpi_1'= \xi_1+\xi_2,\\
\vpi_2'=2\xi_1,\\
\vpi_3'= 2\xi_1+\xi_3+ \xi_4,\\
\vpi_4'= 2(\xi_1+\xi_3),\\
\ldots\\
\vpi'_{2l-1} = 2(\xi_1+\xi_3+\ldots+\xi_{2l-3})+\xi_{2l-1}+
\xi_{2l},\\
\vpi'_{2l} = \xi_1+\xi_3+\ldots+\xi_{2l-3}+\xi_{2l-1}
\end{array}\right.$$

For $n=2l+1$, ~$l>1$, we have

$$\left\{\begin{array}{l}\xi_1= \veps_1 + \veps_2,\\
\xi_2 =  \veps_1 - \veps_2,\\
\ldots,\\
\xi_{2l-1}= \veps_{2l-1} +  \veps_{2l},\\
\xi_{2l} =  \veps_{2l-1} -
\veps_{2l},\\
\xi_{2l+1} =  \veps_{2l+1}.\end{array}\right.\quad\quad\quad
\left\{\begin{array}{l}\vpi_1'= \xi_1+\xi_2,\\
\vpi_2'=2\xi_1,\\
\vpi_3'= 2\xi_1+\xi_3+ \xi_4,\\
\vpi_4'= 2(\xi_1+\xi_3),\\
\ldots\\
\vpi'_{2l-1} = 2(\xi_1+\xi_3+\ldots+\xi_{2l-3})+\xi_{2l-1}+
\xi_{2l},\\
\vpi'_{2l} = 2(\xi_1+\xi_3+\ldots+\xi_{2l-3}+\xi_{2l-1}),\\
\vpi'_{2l+1} = \xi_1+\xi_3+\ldots+\xi_{2l-3}+\xi_{2l-1} +
\xi_{2l+1}.
\end{array}\right.$$
{\textsc Case} $C_n=\mathrm{sp}(2n,K)$, ~$n\geq 3$. Here $m=n$ and
$\vpi'_i=2\vpi_i$ for any $1\leq i\leq n$. We obtain

$$\left\{\begin{array}{l}\xi_1= 2\veps_1,\\
\xi_2 =  2\veps_2,\\
\ldots,\\
\xi_{n}= 2\veps_{n}.\end{array}\right. \quad\quad\quad
\left\{\begin{array}{l}\vpi_1'= \xi_1,\\
\vpi_2'= \xi_1 +\xi_2,\\
\ldots\\
\vpi_n' = \xi_1 +\xi_2+\ldots+\xi_n.\end{array}\right.$$

{\textsc Case} $D_n=\mathrm{o}(2n,K)$, ~ $n\geq 4$. In the case
$n=2l$, we have $m=n$, ~$\vpi'_i=2\vpi_i$ for any  $1\leq i\leq n$
and
$$\left\{\begin{array}{l}\xi_1= \veps_1 + \veps_2,\\
\xi_2 =  \veps_1 - \veps_2,\\
\ldots,\\
\xi_{2l-1}= \veps_{2l-1} +  \veps_{2l},\\
\xi_{2l} =  \veps_{2l-1} - \veps_{2l}.
\end{array}\right.\quad\quad\quad
\left\{\begin{array}{l}\vpi_1'= \xi_1+\xi_2,\\
\vpi_2'=2\xi_1,\\
\vpi_3'= 2\xi_1+\xi_3+ \xi_4,\\
\vpi_4'= 2(\xi_1+\xi_3),\\
\ldots\\
\vpi_{2l-2}'= 2(\xi_1+\xi_3+\ldots+\xi_{2l-3},\\
\vpi_{2l-1}'= \xi_1+\xi_3+\ldots+\xi_{2l-3}+\xi_{2l-1},\\
\vpi_{2l}'=
\xi_1+\xi_3+\ldots+\xi_{2l-3}+\xi_{2l}.\end{array}\right.$$

For $n=2l+1$ we have $m = n-1$,~ $\vpi'_i=2\vpi_i$ for any $1\leq
i\leq n-2$ and $\vpi'_{n-1}= \vpi _n' = \vpi_{n-1} + \vpi_n$. We
obtain

$$\left\{\begin{array}{l}\xi_1= \veps_1 + \veps_2,\\
\xi_2 =  \veps_1 - \veps_2,\\
\ldots,\\
\xi_{2l-1}= \veps_{2l-1} +  \veps_{2l},\\
\xi_{2l} =  \veps_{2l-1} -
\veps_{2l}.\end{array}\right.\quad\quad\quad
\left\{\begin{array}{l}\vpi_1'= \xi_1+\xi_2,\\
\vpi_2'=2\xi_1,\\
\vpi_3'= 2\xi_1+\xi_3+ \xi_4,\\
\vpi_4'= 2(\xi_1+\xi_3),\\
\ldots\\
\vpi'_{2l-1} = 2(\xi_1+\xi_3+\ldots+\xi_{2l-3})+\xi_{2l-1}+
\xi_{2l},\\
\vpi'_{2l} = \xi_1+\xi_3+\ldots+\xi_{2l-3}+\xi_{2l-1}.
\end{array}\right.$$
 {\textsc
Case} $G_2$. We have $m=n=2$, ~~ $\vpi'_i=2\vpi_i$ (for $i=1,2$),
$$
\left\{\begin{array}{l} \xi_1=3\al_1+2\al_2,\\
\xi_2= \al_1.\end{array}\right.
\quad\quad\quad \left\{\begin{array}{l} \vpi'_1 = \xi_1+\xi_2,\\
\vpi'_2 = 2\xi_1.\end{array}\right.$$
 {\textsc Case} $F_4$. We have
$m=n=4$, ~~ $\vpi'_i=2\vpi_i$ (for $i=\overline{1,4}$) and
$$\left\{\begin{array}{l} \xi_1= 2\al_1 + 3\al_2 + 4\al_3 + 2\al_4,\\
\xi_2= \al_2 + 2\al_3 + 2\al_4,\\
\xi_3= \al_2 + 2\al_3,\\
\xi_4= \al_2.\end{array}\right.\quad\quad\quad
\left\{\begin{array}{l} \vpi'_1=2\xi_1,\\
\vpi'_2=3\xi_1+\xi_2+\xi_3+\xi_4,\\
\vpi'_3=2\xi_1+\xi_2+\xi_3,\\
\vpi'_4=\xi_1+\xi_2.\end{array}\right.$$ {\textsc Case} $E_6$. We
have $m=4$ and
$$\left\{\begin{array}{l}
 \xi_1= \al_1 + 2\al_2 + 2\al_3 + 3\al_4 +2\al_5+\al_6,\\
\xi_2= \al_1 + \al_3 + \al_4 + \al_5 + \al_6,\\
\xi_3= \al_3 + \al_4 + \al_5 ,\\
\xi_4= \al_4.\end{array}\right.\quad
\left\{\begin{array}{l} \vpi'_1= \vpi_6' = \vpi_1 + \vpi_6 = \xi_1+\xi_2,\\
\vpi'_2= 2\vpi_1 = 2\xi_1,\\
\vpi'_3= \vpi_5' = \vpi_3 + \vpi_5 = 2\xi_1 + \xi_2 + \xi_3,\\
\vpi'_4= 2\vpi_4 = 3\xi_1 + \xi_2+  \xi_3 +
\xi_4.\end{array}\right.$$ {\textsc Case} $E_7$. We have $m=n=7$,
~~$\vpi'_i=2\vpi_i$ (for $i=\overline{1,7}$),

$$\left\{\begin{array}{l} \xi_1 = 2\al_1 + 2\al_2 + 3\al_3 + 4\al_4 + 3\al_5 +2\al_6 +\al_7,\\
\xi_2 = \al_2 + \al_3 + 2\al_4 + 2\al_5 + 2\al_6  + 2\al_7,\\
\xi_3 = \al_2 + \al_3 + 2\al_4 + \al_5,\\
\xi_4 = \al_7,\\
\xi_5 = \al_2,\\
\xi_6 = \al_3,\\
\xi_7 = \al_5.\end{array}\right.\quad\quad
 \left\{\begin{array}{l}
\vpi'_1= 2\xi_1,\\
\vpi'_2= 2\xi_1 + \xi_2 + \xi_3 + \xi_5,\\
\vpi'_3= 3\xi_1 + \xi_2 + \xi_3 + \xi_6,\\
\vpi'_4= 2(2\xi_1 + \xi_2 + \xi_3),\\
\vpi'_5= 3\xi_1 + 2\xi_2 + \xi_3 + \xi_7,\\
\vpi'_6= 2(\xi_1 + \xi_2),\\
\vpi'_7= \xi_1 + \xi_2 + \xi_4.\end{array}\right.$$

{\textsc Case} $E_8$. We have $m=n=8$, ~~$\vpi'_i=2\vpi_i$ (for
$i=\overline{1,8}$),

$$\left\{\begin{array}{l} \xi_1 = 2\al_1 + 3\al_2 + 4\al_3 + 6\al_4 + 5\al_5 + 4\al_6 + 3\al_7 + 2\al_8,\\
\xi_2 = 2\al_1 + 2\al_2 + 3\al_3 + 4\al_4 + 3\al_5 + 2\al_6  + \al_7,\\
\xi_3 = \al_2 + \al_3 + 2\al_4 + 2\al_5 + 2\al_6 + \al_7,\\
\xi_4 = \al_2 + \al_3 + 2\al_4 + \al_5,\\
\xi_5 = \al_7,\\
\xi_6 = \al_2,\\
\xi_7 = \al_3,\\
\xi_8 = \al_5;\end{array}\right.$$
$$
\left\{\begin{array}{l}
\vpi'_1= 2(\xi_1 + \xi_2),\\
\vpi'_2= 3\xi_1 + 2\xi_2 + \xi_3 + \xi_4 + \xi_4,\hspace{5cm}\\
\vpi'_3= 4\xi_1 + 3\xi_2 + \xi_3 + \xi_4 + \xi_7,\\
\vpi'_4= 2(3\xi_1 + 2\xi_2 + \xi_3 + \xi_4),\\
\vpi'_5= 5\xi_1 + 3\xi_2 + 2\xi_3 + \xi_4+ \xi_8,\\
\vpi'_6= 2(2\xi_1 + \xi_2 + \xi_3),\\
\vpi'_7= 3\xi_1 + \xi_2 + \xi_3 + \xi_5,\\
\vpi'_8 = 2\xi_1.\end{array}\right.$$

\Lemma\Num \label{kij}.~ {\it For any simple Lie algebra and any
$1\leq i\leq m$, the greatest common divider of the row
$(k_{i1},\ldots,k_{im})$ (see formula (\ref{vpk})) is equal to  1
or 2.}\\
\\
 \Proof~ follows from the above formulas of type (\ref{vpk})
for simple Lie algebras  $A_n - E_8$. $\Box$

Denote $k_{ij}'=k_{ij}$, in the case
$\mbox{НОД}(k_{i1},\ldots,k_{im})=1$, and
$k_{ij}'=\frac{1}{2}k_{ij}$, in the case
$\mbox{НОД}(k_{i1},\ldots,k_{im})=2$. Note that in any case
$k'_{ij}\in\Zb$. Denote

\begin{equation}\label{Q} Q_i = Z_1^{k_{i1}'}\cdots Z_m^{k_{im}'}.\end{equation}

 It follows from the formula
(\ref{pz}) that  either $P_i=Q_i$, or $P_i=Q_i^2$. Since $k_{ij}'\in \Zb$, then $Q_i\in\FC^N$.\\
\\
 \Lemma\Num\label{qfield}.\\{\it
 1) ~$\det(k_{ij}') = \pm 1$;\\
 2) ~ $\FC = K(Q_1,\ldots,Q_m)$.}\\
 \\
\Proof~ of item 2) follows from item 1) and corollary \ref{zfield} .
The statement  of item 1) is checked in each case  $A_n - E_8$
separately.
 $\Box$\\
 \\
 \Lemma\Num\label{qA}. ~~{\it  $Q_1,\ldots,Q_m\in \AC^N$.}\\
 \Proof.  As $P_i\in\AC$, then the statement is obvious in the case $P_i=Q_i$.
 Let  $P_i=Q^2_i$. Since $Q_i\in\FC$, ~$P_i\in\AC$, and the ring
 $\AC$ is integrally closed, then $Q_i\in\AC$. As  $P_i\in \AC^N$, then
 $Q_i\in\AC^N$. $\Box$\\
\\
 \Theorem\Num\label{T1}.~ \\{\it
 1)~ The ring of invariants $\AC^N$ of the coadjoint representation of the group $N$ is the polynomial
 ring of
  $Q_1,\ldots,Q_m$.\\
 2)~ For any nonzero elements  $c_1,\ldots, c_m$ of the field  $K$, the
 set,
 defined in $\nog^*$ by the system of equations  $$ Q_1 = c_1, \ldots,
Q_m = c_m,$$ is a coadjoint orbit (of maximal dimension).}\\
\\
 \Proof.
Item  2) follows from the formula (\ref{Q}) and item  3) of
proposition  \ref{LocZ}.

Let us prove item  1). The polynomials  $Q_1,\ldots,Q_m$ lie in
$\AC^N$ (see lemma
 \ref{qA}) and algebraically independent over the field $K$ (see item
 2 of the lemma  \ref{qfield}). To conclude the proof, it is sufficient to show that  $Q_1,\ldots,Q_m$
 generate the ring  $\AC^N$.

 Let a polynomial $F$ lie in  $\AC^N$ and $F$ be an eigenfunction of weight
  $\la$ for the coadjoint representation  $\ad^*_\hog$ of Cartan
  subalgebra.

  Let us prove that the weight  $\la$ is dominant (i.e.  $\la(H_{\al_i})\geq 0$ for any simple root
   $\al_i$).  Identify  $\AC=K[\nog^*]$ with the symmetric algebra  $S(\nog)$. Then  $F$ is contained in  $S(\nog)$
 and  is a weight function for the adjoint representation  $\ad_\bog$ in  $S(\nog)$. Since
$S(\nog)\subset S(\gog)$, then $F$ is a dominant vector for the
adjoint representation of $\gog$ in  $S(\gog)$. This proves that
$\la$ is a dominant weight.

We denote the  $\hog$-weights of polynomials  $Q_1,\ldots,Q_m$ by
$\eta_1,\ldots,\eta_m$. It follows from definition of  $Q_i$ that
either  $\eta_i = \vpi_i'$, or $\eta_i = \frac{1}{2}\vpi_i'$.
Looking through the expressions of  each $\vpi_i'$ in terms of
$\{\vpi_i\}$ in each case  $A_n - E_8$, we obtain that $\eta_i$
coincides with one of the weights  $\vpi_i$, ~ $2\vpi_i$ or $\vpi_i
+ \vpi_{\phi(i)}$ (the last case occurs in the cases  $A_n$,~ $D_n$
(n is odd) and $E_6$). We get  $\eta_i(H_{\al_i}) = 2^\epsilon$,
where $\epsilon$ equals to  0 or 1, and $\eta_i(H_{\al_j}) = 0 $ for
any $j\ne i$ and $1\leq j\leq m$.

Since  $F$ is contained in $\AC^N$, then $F\in\FC^N$. It follows
from corollary  \ref{prod} and item 1) of the lemma  \ref{qfield}
that
$$ F = Q_1^{s_1}\cdots Q_1^{s_m}$$
for some $s_1,\ldots, s_m\in\Zb$. Hence
$$\la = s_1\eta_1 + \ldots + s_m\eta_m.$$
For any simple root  $\al_i$ we obtain
$$\la(H_{\al_i}) = 2^\epsilon s_i \geq 0. $$ We conclude that
$s_i\in\Zb_+$  and, therefore, $F\in K[Q_1,\ldots,Q_m]$. $\Box$

The theorem \ref{T1} implies \\
 \Cor\Num. {\it The polynomials  $Q_1,\ldots,Q_m$
are irreducible over the field  $K$.}

\section{Invariants of the coadoint representation of Borel subalgebras}

By the Killing form, we identify   $\bog^*$ with the lower Borel
subalgebra $\bog_-$. Recall that  by  $\BAC$ we denote the algebra
$K[\bog^*]$ of regular functions on  $\bog^*=\bog_-$. Respectively,
$\BFC$ is the field $K(\bog^*)$ of rational functions on $\bog^*$.
\\
\\
\Theorem\Num \label{B-inv}. {\it The ring of invariants $\BAC^B$ of
the
coadjoint representation of the group  $B$ consists of  $K$.}\\
 \Proof~  Let us  identify the algebra $\BAC$ with the symmetric algebra  $S(\bog)$. Thus, we identify
 the coadjoint representation in  $\BAC$ with the adjoint representation of  $B$ in $S(\bog)$.

Let $F\in S(\bog)^B$. In particular,   $F$ is an invariant of the
adjoint representation of the Cartan subgroup  $H$. Hence,  $F\in
S(\hog)$. Since  $\ad_{e_\al}F = 0$ for any simple root  $\al$, then
$F\in K$.~ $\Box$

Turn to description of the field of invariants $\BFC^B$. Let, as
above,  $w_0$ be the element of the greatest length in the Weyl group.\\
\\
 \Theorem\Num\label{bw}. {\it Let $\gog$ be a simple Lie algebra with  $w_0 = -\mathrm{id}$, then $\BFC^B = K$.}\\
\\
 \Proof.  In the case  $m=n$, where $n$, as above, is the  rank of  $\gog$. Denote by  $h_1,\ldots, h_n$
the dual basis for $\xi_1, \ldots, \xi_n$ in $\hog$. Localization of
the algebra  $\BAC$ with respect to the denominator system
$Z_1,\ldots, Z_n$ coincides with the algebra
\begin{equation}\label{BA1}\Ab_s\otimes \Ab'_n,\end{equation} where
$\Ab_s$ is the algebra of the proposition  \ref{LocZ}, the algebra
$\Ab'_n$ is also a standard Poisson algebra with generators
$$ p_i = Z_i^{-1} h_i, \quad q_i = Z_i,\quad i\in\overline{1,n}.$$ Finally, we have got  $\BFC^B = K$. ~$\Box$

  As in (\ref{Fo}), to a rational function  $F(g)$ on the group $G$ we correspond the formal series
   \begin{equation}\label{Fxt}
 F(\exp t\xt) = t^k(F_0(\xt) + tF_1(\xt)+ t^2
F_2(\xt) +\ldots ), \end{equation} where $\exp t\xt$ is the formal
exponent, $ k\in \Zb$ and the coefficients  $F_i(\xt)$ are rational
functions on $\bog_-$. If $F\in K[G]$, then $F_i(\xt) \in \BAC$.
\\
\\
 \Prop\Num\label{BN}. {\it If $F\in
K[G]^{N\times N}$, then $F_0(\xt) \in
\BAC^N$.}\\
\\
\Proof. Consider the Zarisky-open subset  $B_-N$ in the group $G$.
Similarly to the proof of proposition \ref{NN} we define the
dressing action of the group $N$ on $B_-$. As in formula (\ref{gs}),
we  can show that for any $g\in N$ and $\tilde{s}\in B_-$ the
following equality holds
\begin{equation}\label{gts}
F(\rho_g(\tilde{s})) =  F(\tilde{s}).
\end{equation}
After this the proof is concluded similarly to proposition
\ref{NN}. $\Box$

Let  $\gog$ be a Lie algebra such that  $w_0\ne -\mathrm{id}$. As we
saw above, in this case, the Lie algebra  $\gog$ is either coincide
with $A_n$, or with $D_n$~ ($n$  is odd), or with $E_6$. For this
Lie algebras  $w_0=-\phi$, where $\phi$ is some nontrivial
automorphism of the system of simple roots. Automorphism $\phi$ acts
by permutations on the system of fundamental weights $\vpi_1,\ldots,
\vpi_n$. We consider  $\phi$ to be a permutation of
$\{1,\ldots,n\}$. Choose some subset $\Aog$ in the  set
$\{1,\ldots,n\}$, satisfying the following two properties:\\
i) ~ $\phi(i)\ne i$ for any  $i\in\Aog$;\\
ii)~ for any $i$ from the system $\{1,\ldots,n\}$ either $i\in\Aog$,
or $\phi(i)\in\Aog$.

Easy to see that $|\Aog| = n-m$. We restrict the spherical function
 $S_i(g)$, $i\in\Aog$, on $B$ and decompose as in the formula
(\ref{Fxt}):
\begin{equation}\label{Sxt} S_i(\exp t\xt) = t^k(S_{i0}(\xt) + tS_{i1}(\xt)+ t^2 S_{i2}(\xt)
+\ldots ).\end{equation}

Since any element  $\xt$ of $\bog_- = \bog^*$ is uniquely expressed
in the form  $\xt = x+y$, ~ $y\in\hog$,~ $x\in\nog_-$, then we
consider  any polynomial  $F$ on  $\bog_-$ to be  a polynomial in
two variables $x$ and $y$. We shall write $F(\xt)=F(x)$ (resp.
$F(\xt)=F(y)$), if $F$ does not depend on $y$ (resp. on $x$).
\\
\\
\Prop\Num\label{S12}. {\it We claim that \\
1)~ The zero term   $S_{i0} (\xt) $ of expansion  (\ref{Sxt})
coincides with the zero term  $S_{i0}(x)$
in expansion (\ref{Sx}) for the same $S_i(g)$.\\
2)~ The first term $S_{i1}(\xt)$ of expansion (\ref{Sxt}) can be
represented in the  form
\begin{equation}\label{S1} S_{i1}(\xt) = L_i(y)S_{i0}(x) +
R_i(x),\end{equation} where $\xt = x+y$, ~ $y\in\hog$,~
$x\in\nog_-$,
~ $L_i = \frac{1}{2}(1+w_0)\vpi_i$, and $R_i(x)$ is some polynomial in $x$.}\\
\\
\Proof.  To simplify notations, we denote by  $\vpi$ the $i$th
fundamental weight and by $S(g)$ the spherical function of $i$th's
fundamental representation.  In these notations,  the action of
element $g\in G$ on vector we sign as  $gv$ (instead of $T_gv$).
Respectively, the action of element  $x\in\gog$ we sign as  $xv$
(instead of $\frak{T}_xv$, where $\frak{T} = d_eT$).

Denote by $\exp^{(k)}X$ the series, that is obtained from the series
$\exp X$ by deleting the first  $k$ additives. That is
$$ \exp^{(k)}X = \frac{1}{k!} X^k + \frac{1}{(k+1)!} X^{k+1} +
\ldots .
$$

 It follows from (\ref{Sxt}) that
\begin{equation}\label{ek} S(\exp(t\xt)) = l_0(\exp(t\xt)v_0) =
l_0(\exp^{(k)}(t\xt)v_0).\end{equation}

Decompose  $\xt=x+y$, where $x\in\nog_-$, ~ $y\in\hog$. By
Campbell-Hausdorf formula ~\cite{ser} we have  $$\exp(t\xt) =
\exp(t\xt) \exp(-ty) \exp(ty) =  \exp(t\xt - ty -\frac{t^2}{2}
[\xt,y] + O(t^3)) \exp(ty) = $$ $$ \exp(tx - \frac{t^2}{2} [x,y] +
O(t^3))(1+ty+O(t^2)).$$ Note that
\begin{equation}
\label{ke} \exp^{(k)}(tx - \frac{t^2}{2} [x,y] + O(t^3)) \exp(ty) =
t^k(M_0(\xt) + tM_1(\xt) + O(t^{2})),\end{equation} where
$$\begin{array}{l}M_0(\xt) = M_0(x) =  \frac{1}{k!}x^k,\\
 M_1(\xt) = \frac{1}{k!}(x^ky-\frac{1}{2}(x^{k-1}[x,y] + x^{k-2}[x,y]x+\ldots+[x,y]x^{k-1}))
+ \frac{1}{(k+1)!}x^{k+1} =\\ \frac{1}{2k!} (x^ky+yx^k) +
\frac{1}{(k+1)!}x^{k+1}.\end{array}$$

Taking into account  (\ref{ek}) and (\ref{ke}), we get
\begin{equation}\label{end} S(\exp(t\xt)) = t^k\left(l_0(M_0(x)v_0) +
tl_0(M_1(\xt)v_0) + O(t^2)\right).
\end{equation}
From here $S_0(\xt) = l_0(M_0(x)v_0)$ (this proves statement 1) of
our proposition)  and
$$S_1(\xt) = l_0(M_1(\xt)v_0) = l_0\left(\frac{1}{2k!} (x^ky+yx^k)v_0 +
\frac{1}{(k+1)!} x^{k+1}v_0\right).$$

Since $yv_0=\vpi(y)v_0$ and $l_0(yv) = w_0\vpi(y) l_0(v)$ for any
$v$ from the representation space, then $$ S_1(\xt) =
L_\vpi(y)S_0(x) + R(x),$$ где $$ L_\vpi(y) = \frac{1}{2}
(1+w_0)\vpi(y)~~\mbox{and}~~~ R(x) = l_0\left(\frac{1}{(k+1)!}
x^{k+1}v_0\right).$$ $\Box$

 Denote by  $J_i$  the following rational function

 \begin{equation}\label{Ji} J_i = \frac{S_{i1}}{S_{i0}},~~ \mbox{where}~~
 i\in\Aog.\end{equation}
\\
 \Lemma\Num\label{in}. {\it If $\phi(i)\ne i$, then $J_i$ is an invariant of the coadjoint representation of the group
  $B$.}\\
 \\
 \Proof.~  Any polynomial  $S_{ij}(\xt)$ in decomposition (\ref{Sxt})
 is a eigenvector  for the coadjoint representation of the Cartan subgroup $\hog$ of weight $\vpi'_i$. From here $\ad_\hog^*(J_i)=0$ and, therefore,
 $J_i$ is an invariant for $\Ad^*_H$.

Show that  $J_i$ is an invariant for  $\Ad_N^*$. Really,
 $S_{i0}$ is an invariant by means of  proposition \ref{BN}.

Let us prove that  $S_{i1}$ is also invariant for $\Ad_N^*$. The
fundamental representation  $\vpi_{\phi(i)}$ is conjugate for
$\vpi_i$. Hence, $S_{\phi(i)}(g) = S_i(g^{-1})$. From here
$$ S_i(\exp(t\xt)) = t^k( S_{i0}(\xt) + tS_{i0}(\xt) + O(t^2)),
$$
$$
S_{\phi(i)}(\exp(t\xt)) = S_i(\exp(-t\xt))  = t^k( (-1)^kS_{i0}(\xt)
+ (-1)^{k+1} tS_{i1}(\xt) + O(t^2)).$$

Consider the function
$$ F(g) = \frac{1}{2}(S_i(g) - (-1)^kS_{\phi(i)}(g)).
$$

Expansion (\ref{Fxt}) for the  polynomial  $F(g)$ has the form
$$ F(\exp(t\xt)) = t^{k+1}( S_{i1}(\xt) + O(t)).$$
The polynomial  $F(g)$ is a $N\times N$-invariant; hence, $
S_{i1}(\xt)$ is an invariant for  $\Ad^*_N$ (see proposition
\ref{BN}).~
$\Box$\\
\\
 \Theorem\Num\label{invB}. {\it
If  $\gog$ is a simple Lie algebra,  satisfying  $w_0 \ne
-\mathrm{id}$, then  $\BFC^B$ is the field of rational functions of
the system
$\{J_i(\xt):~ i\in\Aog\}$.}\\
\\
\Proof. First, let us prove that  $w_0(\xi) = -\xi$ for any
$\xi\in\Xi$. Indeed, since  $\xi_1$ is the greatest positive root,
then $w_0(\xi_1) = -\xi_1$. If $\gamma$ is a singular root for
$\xi_1$, then $w_0(\gamma) = -\gamma'$, where $\gamma'$ is also a
singular root for $\xi_1$. Therefore,  $w_0(\De_1^+) = -\De_1^+$. As
a corollary,~ $w_0(\xi_2) = -\xi_2$. Continuing the process further,
we have got $w_0(\xi) = -\xi$ for any $\xi\in\Xi$.

From here $(1+w_0)\Xi=0$ and, therefore,  $(1+w_0)h \perp \Xi$ for
any  $h\in \hog$. Since ~ $L_i= \frac{1}{2}(1+w_0)\vpi_i$ (see item
2) of proposition \ref{S12}), then  $L_i\perp \Xi$. The system
$\{L_i:~ i\in\Aog\}$ for a basis in the orthogonal complement to
$\Xi$ in $\hog$.

By the formulas (\ref{S1}) and (\ref{Ji}), the system  $\{ J_i:~~
i\in \Aog\}$ is algebraically independent.

We complete  $\{L_i:~ i\in\Aog\}$ by the system $h_1,\ldots, h_m$ to
the basis of $\hog$ so that  $\xi_i(h_j) = \delta_{ij}$.

The localization of the algebra  $\BAC$ by the denominator system
$Z_1,\ldots, Z_n$ coincides with the algebra
\begin{equation}\label{BA2}\Ab_s\otimes \Ab'_n\otimes
K[J_i:~~i\in\Aog]\end{equation}
 where $\Ab_s$
is an algebra of proposition  \ref{LocZ}, the algebra  $\Ab'_n$ is
also a standard Poisson algebra with generators
$$ p_i = Z_i^{-1} h_i, \quad q_i = Z_i,\quad i\in\overline{1,n}.$$ In follows that $\BFC^B =
K[J_i:~~i\in\Aog]$. ~$\Box$\\
\\
\Theorem\Num\label{Final}. {\it \\
1) ~Let  $ w_0= -\mathrm{id}$. Then the set, defined in $\bog^*$ by
the system of inequalities   $$Q_1\ne 0,\ldots, Q_m\ne 0,$$
is a coadjoint orbit (of maximal dimension)   of the group $B$ in  $\bog^*$.\\
2)~ Let   $ w_0\ne-\mathrm{id}$.  Then for any system of constants
$\{ c_i\in K:~ i\in\Aog\}$ the set, defined in $\bog^*$ by the
system of relations
$$Q_1\ne 0,\ldots, Q_m\ne 0,$$
$$J_i = c_i,\quad i\in\Aog,$$
is a coadjoint orbit (of maximal dimension)   of the group $B$ in  $\bog^*$.}\\
\\
\Proof~ follows from the presentation of localization of the algebra
 $\BAC$ with respect to $Z_1,\ldots, Z_m$ as the tensor product of type (\ref{BA1}), in the case of item 1), and  of type (\ref{BA2}),
 in the case of item 2).
$\Box$


\begin{thebibliography}{100}


\bibitem{K1}
A.A.Kirillov, Unitary representations of nilpotent Lie groups,
Russian Math.Surveys 17(1962) 53-114.
\bibitem{K-Orb}
A.A.Kirillov, Lecture on the orbit method, Graduate Studiesin Math.,
 64(2004), Rrovidence, RI: AMS.
\bibitem{Deift}
P.Deift, L.C.Li, T.Nanda, C.Tomei, The Toda flow on a generic orbit
is integrable, Comm.Pure Appl. Math. 39(1986) 183-232
\bibitem{Arch}
A.A.Arhangel'skii, Completely integrable hamiltonian systems on
group of triangular matrices,  Sbornik:Mathematics 36(1980) 127-134.
\bibitem{P} A.M.Perelomov, Integrable systems of classical mechaniks and Lie
algebras, Nauka, Moscow, 1990.[in russian]

\bibitem{Trof}
V.V.Trofimov, Euler equation on Borel subalgebras of semisimple Lie
algebras,  Izvestiya:Mathemetics  14(1980)  653-670.

\bibitem{Gekht}
M.I.Gekhtman, M.Z.Shapiro, Noncommutative and Commutative
Integrability of Generic Toda Flows in Simple Lie algebras,  Comm.
Pure Appl. Math  52(1999) 53-84.

\bibitem{Dix}
J.Dixmier, Algebras enveloppanties, Gauthier-Villars, Paris, 1974.
\bibitem{Bur}
N.Burbaki, Groupes et algebras de Lie (chapitre IV-VI), Hermann,
Paris, 1968.
\bibitem{ser}
J.-P. Serre, Lie algebras and Lie groups, Benjamin, New York -
Amsterdam, 1966.







\end{thebibliography}
\end{document}